\documentclass[twocolumn,noversion,nonote,squeezed]{cdcarticle}
\pdfoutput=1
\usepackage{cite}                       
\usepackage{amsmath,amssymb,amsfonts}   
\usepackage[utf8]{inputenc}				      
\usepackage[english]{babel}				      
\usepackage{graphicx}                   
\usepackage{xcolor}
\usepackage{figures/block-diagram}
\usepackage{colonequals}                
\usepackage{float}
\usepackage{enumitem}

\newcommand\Tstrut{\rule{0pt}{2.6ex}}         

\setlist[itemize,enumerate]{noitemsep,nolistsep}

\DeclareTextSymbolDefault{\textbullet}{OMS}

\renewcommand{\epsilon}{\varepsilon}
\newcommand{\real}{\mathbb{R}}
\newcommand{\G}{\mathbf{G}}
\renewcommand{\H}{\mathbf{H}}

\newcommand{\tp}{\mathsf{T}}
\newcommand{\bmat}[1]{\begin{bmatrix}#1\end{bmatrix}}
\newcommand{\grad}{\nabla\!}

\def\neg{\hphantom{-}}

\title{\LARGE\bf
A Tutorial on the Structure of\\Distributed Optimization Algorithms}

\author{Bryan Van Scoy%
\thanks{Department of Electrical and Computer Engineering, Miami University, OH~45056, USA. Email \texttt{bvanscoy@miamioh.edu}}%
\and%
Laurent Lessard%
\thanks{Dept. of Mechanical and Industrial Eng., Northeastern University, MA~02115, USA. Email \texttt{l.lessard@northeastern.edu}\newline
The work of L.~Lessard was supported by the National Science Foundation under Grants No. 2136945 and 2139482.}
}

\note{62\textsuperscript{nd} IEEE Conference on Decision and Control}

\begin{document}

\maketitle

\begin{abstract}
  We consider the distributed optimization problem for a multi-agent system. Here, multiple agents cooperatively optimize an objective by sharing information through a communication network and performing computations. In this tutorial, we provide an overview of the problem, describe the structure of its algorithms, and use simulations to illustrate some algorithmic properties based on this structure.
\end{abstract}

\section{Introduction}

Consider a group of agents that are connected together in a communication network, where each agent is capable of communicating with other agents using the network and performing local computations. For instance, each agent may be a computing node, robot, or mobile sensor.

As an illustrative example, consider the problem of large-scale machine learning. Here, each agent is a computing unit with access to a set of data, and the agents seek to cooperatively build a global model that fits all of the data~\cite{forero2010consensus,johansson2008distributed}. Let $n$ denote the number of agents, and let $f_i$ and $y_i$ denote the loss function and model parameters associated with agent $i\in\{1,\ldots,n\}$. To construct a cohesive global model, we can minimize the total loss over all agents subject to the agents agreeing on the model. This can be formulated as the optimization problem
\begin{subequations}\label{eq:problem}
\begin{alignat}{2}
  &\text{minimize}   \quad && \sum_{i=1}^n f_i(y_i) \label{eq:problem1} \\
  &\text{subject to} \quad && y_1 = y_2 = \ldots = y_n. \label{eq:problem2}
\end{alignat}
\end{subequations}
One approach to solve this problem is for all agents to send their data to a central server and have the server solve the problem to build the model. Some issues with this centralized approach are that i) the computations on the server scale with the number of agents, ii) the system is fragile in that failure of the central server causes the entire system to fail, and iii) data must be transmitted directly over the network and is therefore not private.

\begin{figure}[t]
  \centering\includegraphics{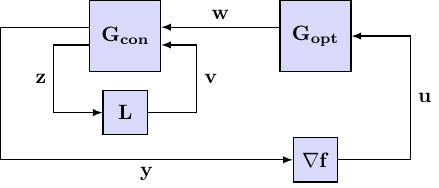}
  \caption{Decomposition of a distributed algorithm into an optimization method $\GOPT$ and second-order consensus estimator $\GCON$.}
  \label{fig:decomp}
\end{figure}

Instead, we seek a \textit{distributed} solution to this problem in which each agent $i$ updates its model $y_i$ using its own loss function $f_i$ and variables that are communicated with neighboring agents. Such algorithms can be made scalable to a large number of agents, robust to failures of individual agents~\cite{self-healing-first-order-distributed-optimization}, and private from other agents~\cite{private-and-hot-pluggable-distributed-averaging}.

Beyond large-scale machine learning, the distributed optimization problem has many other applications, such as multi-agent formation control~\cite{multi-agent-formation-control}, distributed spectrum sensing~\cite{dist_spectrum_sensing}, and distributed allocation of resources~\cite{distributed-resource-allocation-schemes,dist_power_control}.

In this tutorial, our main objectives are as follows:
\begin{enumerate}
  \item Provide an overview of the distributed optimization problem for a multi-agent system.
  \item Describe the structure of algorithms to solve this problem including how they decompose into optimization and consensus components as in Figure~\ref{fig:decomp}.
  \item Use simulations to illustrate some algorithmic properties based on this structural decomposition.
\end{enumerate}
We setup the basic structure of distributed optimization algorithms in Section~\ref{sec:setup} and show how they decompose into optimization and consensus components in Section~\ref{sec:structure}. We then use simulations to illustrate convergence properties in terms of this decomposition in Section~\ref{sec:simulations}.

Throughout the paper, subscripts index the agent and superscripts denote the iteration; for instance, $y_i^k$ is the variable $y$ on agent $i$ at iteration $k$. For a linear time-invariant system $G$, we denote its transfer function as $\widehat G(z)$. We use bold symbols to denote quantities that are aggregated over all agents, such as
\begin{equation}\label{eq:concatenated}
  \y = \bmat{y_1 \\ \vdots \\ y_n} \qquad\text{and}\qquad
  \G = \bmat{G_1 & & \\ & \ddots & \\ & & G_n}.
\end{equation}

\section{Problem setup}\label{sec:setup}

In this section, we describe distributed algorithms for~\eqref{eq:problem} in which agents compute their local gradient and share information through a communication network.

\subsection{Communication network}

\begin{figure}[H]
  \centering\includegraphics{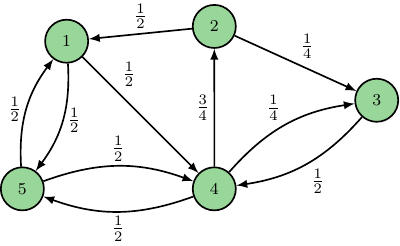}
  \caption{A weighted directed graph that represents the communication network among a group of agents.}
  \label{fig:graph}
\end{figure}

We describe the communication network among the agents as a weighted directed graph, such as in Figure~\ref{fig:graph}. Each vertex in the graph corresponds to an agent in the network and is represented by a circle. Edges in the graph are represented by arrows and indicate the flow of information from one agent to another. The weight of an edge is the amount by which the information is weighted.

Suppose that each agent $i$ has a scalar variable $z_i$. One way for an agent to fuse its information with that of its neighbors is to compute a weighted average of the difference between its variable and those of its neighbors. The weight $a_{ij}$ that agent~$i$ places on information from agent~$j$ is the weight of the edge from node~$i$ to node~$j$ in the graph. This is a linear operation over the concatenated vector $\z = (z_1,\ldots,z_n)\in\real^n$ that is represented by multiplication with the Laplacian matrix $L\in\real^{n\times n}$. In particular, the $i\textsuperscript{th}$ component of the product is
\begin{equation}\label{eq:L}
  (L\z)_i = \sum_{j=1}^n a_{ij}\,(z_i-z_j).
\end{equation}
The weight $a_{ij}$ is nonzero only when agent $j$ is able to send information to agent $i$, so agent $i$ can compute this quantity using the communication network. For the graph in Figure~\ref{fig:graph}, the Laplacian matrix is
\[
  L = \bmat{
      \neg 1 &     -\tfrac{1}{2} & \neg 0 & \neg 0 &     -\tfrac{1}{2} \\[1mm]
      \neg 0 & \neg \tfrac{3}{4} & \neg 0 &     -\tfrac{3}{4} & \neg 0 \\[1mm]
      \neg 0 &     -\tfrac{1}{4} & \neg \tfrac{1}{2} &     -\tfrac{1}{4} & \neg 0 \\[1mm]
          -\tfrac{1}{2} & \neg 0 &     -\tfrac{1}{2} & \neg \tfrac{3}{2} &     -\tfrac{1}{2} \\[1mm]
          -\tfrac{1}{2} & \neg 0 & \neg 0 &     -\tfrac{1}{2} & \neg 1}.
\]
This graph is balanced in that, for each node, the sum of the weights of all incoming edges is equal to that of the outgoing edges~\cite{plp}. In terms of the Laplacian matrix, this means that each row and column sums to one. Balanced graphs preserve averages since $\sum_{i=1}^n (L\z)_i = \sum_{i=1}^n z_i$.

\begin{figure}[htb]
  \centering\includegraphics{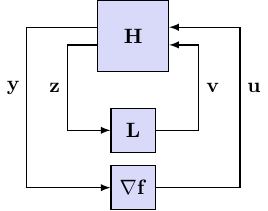}
  \caption{Structure of a general distributed algorithm.}
  \label{fig:distralg}
\end{figure}

\subsection{Distributed algorithms}\label{sec:distralg}

\begin{subequations}\label{eq:distralg}
We now describe the structure of distributed algorithms. Each agent $i$ maintains an estimate $y_i$ of the optimal solution to~\eqref{eq:problem} and can evaluate its local gradient to obtain the quantity
\begin{equation}\label{eq:distralg-gradient}
  u_i = \grad f_i(y_i).
\end{equation}
Agent $i$ can also communicate some quantity $z_i$ with its local neighbors and fuse the information using the graph Laplacian to obtain
\begin{equation}\label{eq:distralg-communication}
  v_i = \sum_{j=1}^n a_{ij}\,(z_i-z_j).
\end{equation}
And finally, the algorithm must determine how to choose the point $y_i$ at which to evaluate the gradient and the point $z_i$ to communicate with neighboring agents in the network. We assume each agent uses the same algorithm and represent this operation as
\begin{equation}\label{eq:distralg-alg}
  \bmat{y_i \\ z_i} = H \bmat{u_i \\ v_i}.
\end{equation}
We focus on algorithms for which $H$ is causal and linear time-invariant (LTI), although some algorithms in the literature are nonlinear and/or time-varying~\cite{distributed-nesterov}. Together, equations \eqref{eq:distralg-gradient}--\eqref{eq:distralg-alg} represent the algorithm on agent $i$.
\end{subequations}

\begin{subequations}
We can represent distributed algorithms compactly using the concatenated vectors~$\u$,~$\v$,~$\y$, and~$\z$ as in~\eqref{eq:concatenated}. In terms of these concatenated vectors, \eqref{eq:distralg-gradient} becomes
\begin{equation}
  \u = \df(\y) \quad\text{where}\quad \df = \text{diag}(\grad f_1,\dots,\grad f_n).
\end{equation}
Likewise, using the Kronecker product $\otimes$, equation \eqref{eq:distralg-communication} can be represented as
\begin{equation}
  \v = \L\z \quad\text{where}\quad \L = L\otimes I_m
\end{equation}
with $m$ the dimension of the communicated variable $z_i$. And finally, \eqref{eq:distralg-alg} becomes
\begin{equation}
  \bmat{\y \\ \z} = \H \bmat{\u \\ \v} \ \text{where}\ \H = \bmat{I_n\otimes H^{11} & I_n\otimes H^{12} \\ I_n\otimes H^{21} & I_n\otimes H^{22}}.
\end{equation}
\end{subequations}
These relationships are summarized by the block diagram in Figure~\ref{fig:distralg}, where the system $\H$ is in feedback with the gradient $\df$ and Laplacian $\L$.

\section{Algorithm structure}\label{sec:structure}

The previous section describes the basic structure of a distributed algorithm as a system $H$ in feedback with the gradient of the objective functions and the Laplacian matrix. It does not, however, provide any insight into how properties of the algorithm depend on~$H$, or even what choices of $H$ lead to sensible algorithms\footnote{For instance, a desireable property is for all fixed points of the algorithm to correspond to optimal solutions of~\eqref{eq:problem}.}.

Intuitively, the optimization problem~\eqref{eq:problem} is a combination of optimization (minimizing the sum of the functions) and consensus (having the agents agree on the solution). In this section, we first review algorithms for optimization and consensus separately, and then describe how any distributed algorithm of the form in Sec.~\ref{sec:distralg} decomposes into these two components.

\subsection{Consensus estimators}\label{sec:con}

We now describe the problem of \textit{consensus}. Suppose each agent $i$ observes a (potentially time-varying) signal $w_i^k$. A consensus estimator is an iterative procedure for each agent to estimate the average signal $w_\text{avg}^k = \frac{1}{n}\sum_{i=1}^n w_i^k$ by sharing information with its local neighbors~\cite{consensus}.

\begin{figure}[htb]
  \centering\includegraphics{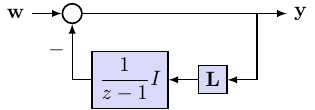}
  \caption{The proportional (P) estimator.}
  \label{fig:P-estimator}
\end{figure}

The block diagram of one particular consensus estimator is shown in Figure~\ref{fig:P-estimator}.

Let $x_i$ denote the state of the estimator on agent $i$. Then the proportional estimator is described by the recursion
\begin{subequations}\label{eq:P-estimator}
\begin{align}
  y_i^k &= w_i^k - x_i^k \label{eq:P-estimator2} \\
  x_i^{k+1} &= x_i^k + \sum_{j=1}^n a_{ij}\,\bigl(y_i^k - y_j^k\bigr) \label{eq:P-estimator1}
\end{align}
\end{subequations}
where the state is initialized such that $\sum_{i=1}^n x_i^0 = 0$.

\begin{figure}[htb]
  \begin{minipage}{0.55\linewidth}
    \centering\includegraphics{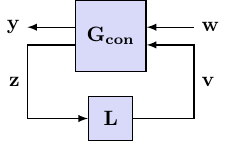}
  \end{minipage}%
  \begin{minipage}{0.45\linewidth}
    \begin{align*}
      \bmat{y_i \\ z_i} &= \Gcon \bmat{w_i \\ v_i} \\[5pt]
      v_i &= \sum_{j=1}^n a_{ij}\,(z_i-z_j)
    \end{align*}
  \end{minipage}
  \caption{General form of a consensus estimator.}
  \label{fig:consensus}
\end{figure}

In general, each agent $i$ has an input signal $w_i$ for which the agents seek to compute the average, an output signal $y_i$ that estimates the average of the inputs, a signal $z_i$ that the agent communicates with neighbors, and a signal $v_i$ that is the result of applying the Laplacian matrix to the communicated variables. The block diagram of a general consensus estimator with these components is shown in Figure~\ref{fig:consensus}.

For example, the transfer function of the P estimator is
\begin{equation}
  \Gconz = \bmat{1 & \frac{-1}{z-1} \\ 1 & \frac{-1}{z-1} \Tstrut}.
\end{equation}
When the input signal is constant, the output of the P estimator converges asymptotically to the average of the input signal, that is, it has zero steady-state error~\cite{consensus}. Such estimators are called \textit{first-order} estimators. Likewise, a \emph{second-order} estimator asymptotically tracks the average of signals whose deviations from their average are ramps. One way to construct a second-order estimator is by combining two first-order estimators in series.

\subsection{Optimization methods}\label{sec:opt}

Now consider a single agent $i$ that seeks to optimize its own objective function $f_i$ (as opposed to the sum of all the functions). To do so, the agent may use a gradient-based optimization method~\cite{lessard} that sequentially queries its gradient $\grad f_i$.

\begin{figure}[htb]
  \begin{minipage}{0.6\linewidth}
    \centering\includegraphics{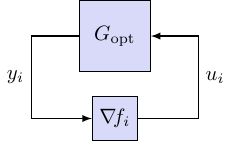}
  \end{minipage}%
  \begin{minipage}{0.4\linewidth}
    \begin{align*}
      y_i &= \Gopt\,u_i \\[5pt]
      u_i &= \grad f_i(y_i)
    \end{align*}
  \end{minipage}
  \caption{General form of an optimization method.}
  \label{fig:optimization}
\end{figure}

A particular optimization method is the \textit{gradient method}, which is described by the recursion
\begin{equation}\label{eq:GM}
  y_i^{k+1} = y_i^k - \alpha\,\grad f_i(y_i^k)
\end{equation}
where $\alpha>0$ is the stepsize. In general, an optimization method applies a discrete-time dynamical system $\Gopt$ to the signal of gradient values $u_i$ to choose the point $y_i$ at which to evaluate the next gradient. The block diagram of a general optimization method is shown in Figure~\ref{fig:optimization}.

For example, the transfer function of the gradient method is
\[
  \Goptz = \frac{-\alpha}{z-1}.
\]
For the optimization method to have a fixed point that satisfies the first-order optimality conditions (that is, the gradient is zero), the transfer function must have a pole at $z=1$. The gradient method uses the minimal number of states to satisfy this requirement, but there are also \textit{accelerated} methods that use additional states to achieve faster convergence~\cite{tmm}.

\subsection{Structural decomposition}\label{sec:decomp}

Since optimization~\cite{lessard} and consensus~\cite{consensus} have been well-studied in the literature, we seek to represent distributed algorithms as a combination of these two components. It turns out that this is possible: we can combine a valid optimization method and second-order consensus estimator using Figure~\ref{fig:decomp} to form a valid distributed algorithm~\cite[Theorem 2]{optcon}. Conversely, any valid distributed algorithm decomposes in this way~\cite[Theorem 1]{optcon}.

While quite general, these results hold for algorithms in which $\Gopt$, $\Gcon$, and $H$ are causal LTI systems that satisfy certain properties to be \textit{valid}. The notion of a valid algorithm varies depending on the type, but roughly means that the algorithm behaves as desired in the simplest scenario. A valid optimization method, for instance, must converge to the optimal solution when applied to the quadratic function $y\mapsto\tfrac{\epsilon}{2} \|y-y^\star\|^2$ for all $y^\star$ and all $\epsilon>0$ sufficiently small; see~\cite{optcon} for additional details.

\section{Simulations}\label{sec:simulations}

We now use simulations to illustrate some properties of distributed algorithms based on the decomposition in Figure~\ref{fig:decomp}. We first setup the problem and then describe each of the various properties.

Consider a set of $n=5$ agents that are connected in a communication network as shown in Figure~\ref{fig:graph}. Suppose the agents cooperate to solve a machine learning problem in which the loss function is quadratic in the model parameters. That is, the agents solve a distributed linear least squares problem. The objective function on agent $i$ is the quadratic
\[
  f_i(y) = \tfrac{1}{2} y^\tp A_i y - b_i^\tp y
\]
parameterized by the symmetric matrix ${A_i\in\real^{d\times d}}$ and the vector $b_i\in\real^d$. The gradient is the linear function
\[
  \grad f_i(y) = A_i y - b_i.
\]
To generate the data, we sample the matrix $A_i$ such that its eigenvalues are evenly spaced in the interval $[\tfrac{1}{10},1]$, and we sample each element of the vector $b_i$ from a standard normal distribution.

The dimension of the model parameters is $d=3$. Since the results depends on the problem data which is random, we simulate $1000$ trials for each scenario.

\subsection{Optimization and consensus errors}

At each iteration, the error is a measure of the distance between the iterates of all the agents and the optimal solution to the distributed optimization problem~\eqref{eq:problem}. To define the error, we use the first-order optimality conditions which are as follows:

\begin{subequations}\label{eq:optimality}
\begin{itemize}
  \item Optimality: the sum of the gradients is zero
  \begin{equation}\label{eq:optimality-opt}
    \sum_{i=1}^n \grad f_i(y_i) = 0
  \end{equation}
  \item \textbf{Consensus:} the agents agree on the optimizer
  \begin{equation}\label{eq:optimality-con}
    y_1 = y_2 = \ldots = y_n
  \end{equation}
\end{itemize}
\end{subequations}

We characterize the error of the iterates in terms of their distance from satisfying these conditions. This consists of two components: the size of the average gradient (the optimization error) and the amount of disagreement (the consensus error). For iterates $y_1^k,\ldots,y_n^k$, we define the optimization and consensus errors as follows:
\[
  e_\text{opt}^k = \biggl\|\sum_{i=1}^n \grad f(y_i^k)\biggr\|
  \quad\text{and}\quad
  e_\text{con}^k = \sum_{i=1}^n\,\biggl\|y_i^k - \frac{1}{n}\sum_{j=1}^n y_j^k\biggr\|.
\]
We take the total error as the maximum of the optimization and consensus errors, $e^k = \max\{e_\text{opt}^k,e_\text{con}^k\}$, which is zero if and only if the first-order optimality conditions are satisfied.

\begin{figure}[htb]
  \centering\includegraphics[width=\linewidth]{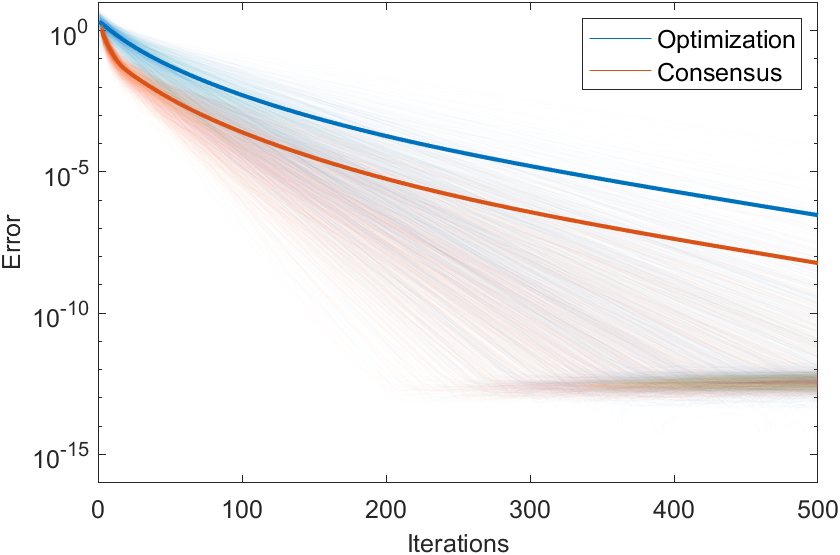}
  \caption{Optimization and consensus errors.}
  \label{fig:opt-con}
\end{figure}

Figure~\ref{fig:opt-con} shows the optimization and consensus errors for each trial (thin) as well as the mean (thick) as a function of the number of iterations. Here, we use the distributed optimization algorithm in Figure~\ref{fig:decomp} in which $\Gcon$ is two P estimators connected in series and $\Gopt$ is the gradient method with stepsize~$\alpha = 0.25$.

\subsection{Factored form}

\begin{figure}[htb]
  \centering\includegraphics[width=\linewidth]{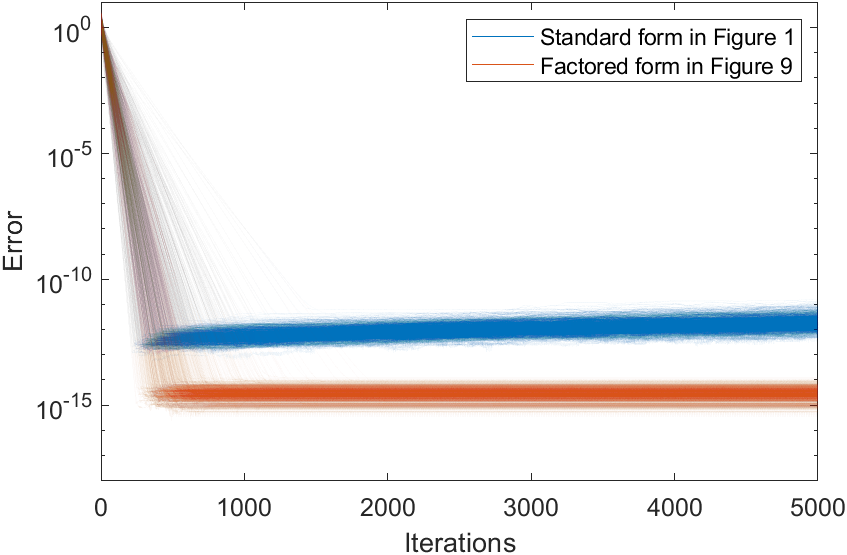}
  \caption{The factored form is numerically stable while the general form is not.}
  \label{fig:factored}
\end{figure}

Suppose that we iterate the algorithm for a long time. From the first-order optimality condition~\eqref{eq:optimality-opt}, the average gradient must be zero at the optimal solution. The gradient of each individual agent, however, is nonzero in general\footnote{If the gradient of each agent were zero at the optimal solution, then there would be no need to cooperate since each agent could solve the global problem by minimizing its local function!}. Recall that the optimization method must have a pole at $z=1$. The nonzero constant gradient resonates with this pole which causes $w_i$ to grow as a ramp. As this signal grows without bound, the error also grow over time as shown in Figure~\ref{fig:factored} (blue).

\begin{figure}[htb]
  \centering\includegraphics{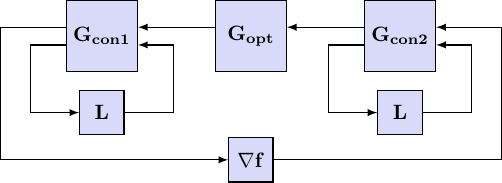}
  \caption{Factored form, where the second-order consensus estimator factors into two first-order estimators.}
  \label{fig:factored-diagram}
\end{figure}

We can fix this issue, however, if the consensus estimator \textit{factors} into two first-order estimators as shown in Fig.~\ref{fig:factored-diagram}.

In this case, the nonzero gradient on each agent is first averaged by the estimator before being applied to the optimization method. Since the average gradient \textit{is} zero, the input to the optimization method is now zero at the optimizer, so the signals no longer grow over time. This is illustrated in Figure~\ref{fig:factored} (red), where the error remains at the numerical precision of the computer over time. In this case, both consensus estimators are the P estimator.

\subsection{Accelerated convergence}

\begin{figure}[htb]
  \centering\includegraphics{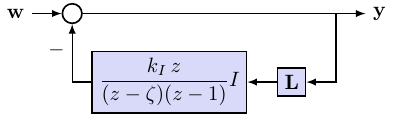}
  \caption{Accelerated consensus estimator.}
  \label{fig:P-estimator-accelerated}
\end{figure}

The decomposition in Figure~\ref{fig:decomp} provides an intuitive procedure to accelerate the convergence of the algorithm. To accelerate the convergence, we can replace the consensus estimator with the accelerated version in Figure~\ref{fig:P-estimator-accelerated} that uses additional dynamics to (potentially) accelerate the rate of convergence.

\begin{figure}[htb]
  \centering\includegraphics[width=\linewidth]{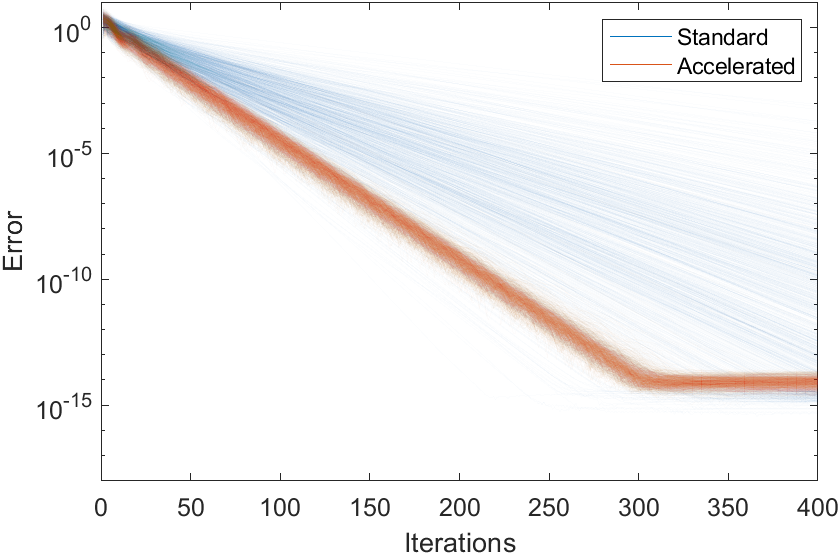}
  \caption{Combining an accelerated consensus estimator and accelerated optimization method may lead to faster convergence.}
  \label{fig:accelerated}
\end{figure}

Similarly, we can replace the gradient method with an accelerated optimization method. Many common first-order methods have the transfer function
\[
  \Goptz = \frac{-\alpha\,(z+\gamma\,(z-1))}{(z-\beta)(z-1)}.
\]
As our intuition suggests, using the accelerated consensus estimator and optimization method improves the convergence rate as shown in Figure~\ref{fig:accelerated} (red), where we use parameters $(\alpha,\beta,\eta) = (0.1,0.8,0)$ and $(\zeta,k_I) = (0.1,1.1)$.

\subsection{Robustness}

\begin{figure}[H]
  \centering\includegraphics{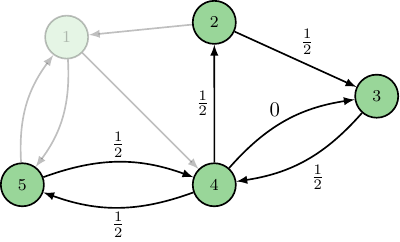}
  \caption{The modified graph after agent $1$ leaves the communication network.}
  \label{fig:graph-drop}
\end{figure}

Now suppose agent $1$ leaves the network; the agent may have malfunctioned, ran out of power, or been hijacked by an adversary. The modified graph is shown in Fig.~\ref{fig:graph-drop}, where agent $1$ and all of its connected edges are opaque to symbolize that it no longer affects the computation\footnote{To maintain a balanced graph, the other agents update their weights so that the sum of the incoming weights is equal to that of the outgoing weights. While agent~$4$ is still capable of sending information to agent~$3$, the corresponding weight is zero indicating that the information is unused.}.

\begin{figure}[htb]
  \centering\includegraphics[width=\linewidth]{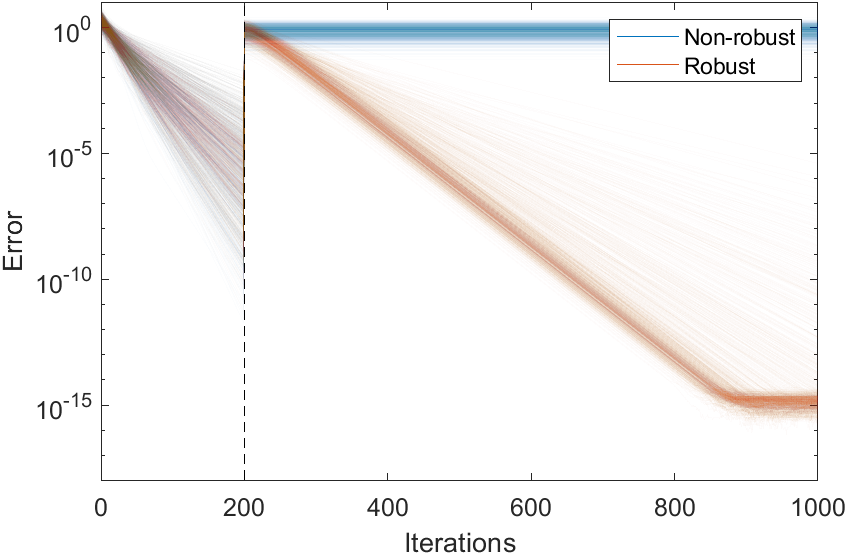}
  \caption{Agent $1$ leaves the communication network at iteration $k=200$. Using a robust consensus estimator enables the algorithm to recover from such changes.}
  \label{fig:robust}
\end{figure}

The proportional consensus estimator requires specific initialization in that the average state must be zero. This average is invariant in that it does not change over time. From~\eqref{eq:P-estimator2}, we observe that the average state $x_i$ must be zero for the average output $y_i$ to equal the average of~$w_i$. While we initially set the average state to be zero, it is in general nonzero once the network changes. This results in a systemic error as shown in Figure~\ref{fig:robust} (blue).

\begin{figure}[htb]
  \centering\includegraphics{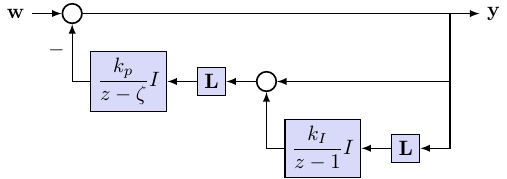}
  \caption{The proportional--integral (PI) estimator, which is robust to changes in the network.}
  \label{fig:PI-estimator}
\end{figure}

An algorithm is \textit{robust} to changes in the communication network if it does not require a specific initialization, in which case it eventually recovers from such changes. To obtain an algorithm that is robust, we can simply replace the consensus estimator $\Gconright$ in Fig.~\ref{fig:factored} with an estimator that is robust. One such estimator is the PI estimator whose block diagram is shown in Figure~\ref{fig:PI-estimator}.

Using this estimator, the error has a transient after the change in the communication network but then converges to zero as shown in Figure~\ref{fig:robust} (red). Here, we use the parameters $(k_p,k_I,\zeta) = (1,0.5,0.95)$.

\section{Conclusion}

In this tutorial, we studied distributed optimization for a multi-agent system. We described the structure of algorithms in terms of optimization methods and consensus estimators, and we illustrated some of the properties through simulations on a machine learning problem.

We focused on algorithms that consist of a linear time-invariant system in feedback with the gradient and the Laplacian; however, other types of algorithms may be needed depending on the class of objective functions to be optimized. We also focused on certain properties that depend on the algorithm structure, but we did not describe how to systematically analyze the convergence rate to the optimal solution; such analysis can be done using the techniques in~\cite{shuo,SVL,sundararajan_allerton}.

\bibliographystyle{IEEEtran}

\setlength{\parskip}{0pt plus 0.3ex}
{\footnotesize\bibliography{references}}

\end{document}